\documentclass[journal,twoside,web]{ieeecolor}
\usepackage{generic}
\usepackage{cite}
\usepackage{amsmath,amssymb,amsfonts}
\usepackage{algorithmic}
\usepackage{graphicx}
\usepackage{algorithm,algorithmic}
\usepackage{hyperref}
\hypersetup{hidelinks=true}
\usepackage{textcomp}

\usepackage{mathdots}
\usepackage{mathbbol}
\usepackage{comment}

\def\BibTeX{{\rm B\kern-.05em{\sc i\kern-.025em b}\kern-.08em
    T\kern-.1667em\lower.7ex\hbox{E}\kern-.125emX}}
\begin{document}
\title{Spectral Decomposition of Discrete-Time Controllability Gramian and Its Inverse via System Eigenvalues} 
\author{Iskakov Alexey 
\thanks{Iskakov A.B. is with the V. A. Trapeznikov Institute of Control Sciences of Russian Academy of Sciences, 
65 Profsoyuznaya street, Moscow 117997, Russia (e-mail: iskalexey@gmail.com).}
}

\maketitle

\begin{abstract}
This paper develops a closed-form spectral decomposition framework for the  Gramian matrices of discrete-time linear dynamical systems. The main results provide explicit decompositions of the discrete-time controllability Gramian and its inverse in terms of the eigenvalues of the dynamics matrix, yielding a mode-resolved representation of these matrices. In contrast to the more common use of aggregate Gramian characteristics, such as eigenvalues, singular values, determinants, and trace-based metrics, the proposed approach describes the internal structure of the Gramian itself through contributions associated with individual modes and their pairwise combinations. 
The framework is extended further to the solution of the discrete-time Lyapunov difference equation, placing the obtained formulas in a broader 
context relevant to the analysis and computation of time-varying and nonlinear systems.
In addition, the decomposition is generalized to systems whose dynamics matrix has multiple eigenvalues, 
enabling a closed-form estimation of the effects of resonant interactions between eigenmodes.
The proposed results provide a structural tool for the analysis of controllability, observability and stability in discrete-time systems and complement existing Gramian-based methods used in model reduction, estimation, actuator and sensor selection, and energy-aware control. Beyond their theoretical interest, the derived decompositions may support the development of improved computational procedures and more informative performance criteria for a range of discrete-time control problems.
\end{abstract}

\section{Introduction}

{\it Controllability and observability Gramians} are fundamental objects in discrete-time system theory because they connect structural properties with quantitative performance. For stable discrete-time LTI systems, these matrices encode least-energy reachability and output sensitivity, providing quantitative metrics for analysis, estimation, and synthesis. In recent work they appear in model reduction, sensing and actuation design, and the energy-aware analysis of large-scale and networked systems \cite{Bag-2009},\cite{Chahl-2012},\cite{Pasq-2013},\cite{Summ-2014}. 
A substantial body of recent literature uses discrete-time Gramians through spectral or scalar quantities derived from them. In model reduction, balanced truncation 
relies on Hankel singular values and the spectral decay of Gramian solutions to identify low-order approximations and to establish error bounds \cite{Chahl-2012}, \cite{Benn-2016}, \cite{Duff-2021}, \cite{Kursch-2024}. In estimation and sensing, observability-type Gramians are tied to Kalman filtering, Fisher information, and sensor selection criteria \cite{Bag-2009}, \cite{Tzoum-2015}, \cite{Boy-2025}. In actuator placement and network control, scalar metrics derived from the controllability Gramian, including trace, determinant, minimum eigenvalue, and inverse-based measures, are widely used to assess control difficulty and to guide design \cite{Summ-2014}, \cite{Cor-2014}, \cite{Sia-2018}, \cite{Bag-2022}. Large-scale numerical methods for Stein and Lyapunov equations have further strengthened the practical role of Gramian techniques in high-dimensional settings \cite{Li-13}, \cite{Sad-12}.

This broad literature uses {\it spectral properties of discrete-time Gramians} through aggregate quantities, such as its eigenvalues, determinant, trace, balancing-related singular values, or related summaries of control energy and estimation quality \cite{Benn-2016}, \cite{Bag-2022}, \cite{Yan-2015}. There are also some papers relating such quantities to spectral features of the dynamics matrix and to energy scaling laws in networked systems \cite{Yan-2012}, \cite{Bag-2021}. However, the literature rarely resolves the controllability or observability Gramian into explicit modal contributions attached directly to the eigenvalues of the dynamics matrix. A few papers come close. \cite{Dilip-2019} derives a Hadamard-product representation of the controllability Gramian for diagonalizable LTI systems and exploits it for actuator and sensor selection, while \cite{Bag-2022} obtains a closed-form determinant expression that exposes how reachable volume depends on spectral properties of the dynamics. 
In \cite{Yad-2016} and \cite{DAN-2017}, decompositions of solutions to the discrete algebraic Lyapunov equation were proposed in terms of the dynamics matrix spectrum. However, these results have not yet led to the development of a widely accepted decomposition framework in which the discrete-time Gramian is divided into interpretable pieces that are directly linked to the eigenstructure of the dynamics matrix. Moreover, the inverse matrix of the Gramian has not been studied from this viewpoint at all. Recent continuous-time studies suggest that such decompositions may carry useful structural information \cite{Yad-2022}, \cite{Isk-2025}, but comparable discrete-time formulations have not yet been developed. This leaves room for new work that complements existing Gramian-based metrics with a finer description of how individual modes and modal interactions shape controllability and observability.

This paper addresses this gap. It derives a closed-form {\it spectral decomposition via system eigenvalues} (SDSE) for the discrete-time Gramian and its inverse, yielding a mode-resolved description of controllability and observability beyond aggregate spectral summaries. By controllability–observability duality, the corresponding results for the observability Gramian follow directly, and hence, are omitted. {\it Controllability and observability canonical forms}, which represent system dynamics using a minimal number of nonzero parameters, provide a convenient framework for control and observer design. Accordingly, we analyze Gramians expressed in controllability or observability canonical form, yielding particularly simple and compact spectral decompositions that are independent of the input and output matrices. The main contributions of this work are summarized as:
\vskip 1mm
\begin{itemize}
\item
Closed-form SDSE of the inverse of the discrete-time controllability Gramian are derived in Theorems~5,~6 and~8, providing a quantitative characterization of the contributions of individual eigenmodes and associated system components to the minimum control energy.
\vskip 1mm
\item 
The existing SDSE results derived in \cite{Yad-2016}, \cite{DAN-2017} have been extended to include solutions to the difference Lyapunov equation with an arbitrary initial condition (Theorems~1 and~3), enabling time-varying and nonlinear cases to be treated, as well as transient analysis to be performed.
\vskip 1mm
\item 
The SDSE of the Gramian and its inverse are generalized to systems with multiple eigenvalues in the dynamics matrix in Theorems~7 and~8, allowing closed-form analysis of resonant interactions between eigenmodes.
\end{itemize}
These results provide an explicit mode-resolved description of discrete-time Gramian matrices that complements the existing literature on Gramian-based metrics, spectral summaries, and energy bounds \cite{Bag-2022}, \cite{Dilip-2019}, \cite{Yad-2016}, \cite{Isk-2025}.

The remainder of the paper is organized as follows. Section II introduces the problem formulation and reviews relevant results on spectral decompositions of Gramians. Section III derives SDSE for the solution to the difference Lyapunov equation with arbitrary initial conditions. Section IV presents SDSE for the solutions of both algebraic and difference Lyapunov equation in companion form. Section V obtains SDSE of the inverse Gramian in controllability canonical form. Section VI extends the analysis to systems with multiple eigenvalues. Section VII presents the conclusions of the paper.\\

\section{Problem Formulation and Preliminaries}

A {\it discrete time LTI dynamical system} is described by the equations
\begin{equation}  \label{sys-0}
x(t+1) =A x(t) + B u(t), \ \ t \in  \mathbb{Z},
\end{equation}
where $x(t)\in \mathbb{R}^n$ is the state of the system, $u(t)\in \mathbb{R}^m$ is the input vector, and the constant matrices $A\in \mathbb{R}^{(n\times n)}$ and  $B\in \mathbb{R}^{(n \times m)}$ describe system dynamics and input, respectively. 

The {\it discrete-time Lyapunov algebraic and difference equations (also known as Stein equations)} for system \eqref{sys-0} are
\begin{equation}  \label{Lyap-algeb} 
 A P A^T - P = - BB^T ,                                      
\end{equation}
\begin{equation}  \label{Lyap-diff} 
P(t+1) = A P(t) A^T + B B^T , \ \ P(0) =P_0  , \ \  t \in  \mathbb{Z}.                                
\end{equation}

Let matrix $A$ have {\it the characteristic equation} $\mathrm{det}(Is - A) = \sum^n_{i=0} a_i s^i$, where $a_n=1$, and let $\mathcal{C} = [ B,AB,A^2 B,...,A^{(n-1)} B]$ be {\it the controllability matrix}. Consider a single-input (SI) controllable system where $m=1$ and $B=b\in R^n$. 
Then, according to \cite{Hauksdottir-2009}, there exists a non-singular similarity matrix 
\begin{equation} \label{transfer-matrix}
T = \mathcal{C} \mathcal{H}_u, \ \ \mathcal{H}_u =      
\begin{bmatrix}
a_1 & a_2 & \cdots & a_{n-1}  & 1  \\ 
a_2 & \iddots & \iddots  & 1            & 0  \\
\vdots &  \iddots   & \iddots & \iddots & \vdots \\
a_{n-1} & 1 & 0    &  \iddots      & 0  \\
1        & 0   &  \cdots      &  0     & 0  
\end{bmatrix} ,
\end{equation}
that transforms \eqref{sys-0} into {\it the controllability canonical form}
\begin{gather} 
x_c (t+1) = A_C \, x_c(t) + b_C \, u(t) , \ \ \text{where} \ \nonumber \\              
x = T x_c, \  A \, T = T A_C, \  B = T b_C, \label{sys-c} \\   
A_C = 
\begin{bmatrix}
0_{(n-1)\times1} &  & I_{(n-1)\times(n-1)} &   \\ 
-a_0 & -a_1 & \dots & -a_{n-1}              
\end{bmatrix} ,  \nonumber \\
b_C= [ 0, 0, \cdots, 0, 1]^T. \nonumber
\end{gather}
Lyapunov equations \eqref{Lyap-algeb} and \eqref{Lyap-diff} are transformed to 
\begin{gather}    
 A_C P_C A^T_C - P_C = - \, b_C b_C^T ,  \label{Lyap-algeb-c} \\ 
P_C(t+1) = A_C P_C(t) A^T_C + \, b_C b_C^T , \ \ P_C(0)= P^C_{0},  \label{Lyap-diff-c} \\  
P = T P_C \, T^T .    \label{P-to-orig-sys}  
\end{gather}

To  find the controllability Gramian $P$ from $P_C$ in the general case of a multiple-input (MI) system, 
each of the $m$ columns of $B$ can be treated separately using the expression proposed in \cite{Hauksdottir-2009}:
\begin{equation}  \label{MI-P-from-Pc}  
P = \mathcal{C} (\mathcal{H}_u P_C \mathcal{H}_u  \otimes I_m) \, \mathcal{C}^* ,    
\end{equation}
where $I_m$ is a unit matrix $m\times m$, $(\cdot)^*$ is the complex conjugate transpose, and $\otimes$ is the matrix Kronecker product. \\

In the subsequent presentation, we use the following result as a starting point. 
The infinite controllability Gramian $P$ of stable system \eqref{sys-0} can be obtained as a solution of the algebraic Lyapunov equation \eqref{Lyap-algeb}
and represented as a sum of Hermitian matrices corresponding to either individual eigenvalues or their pairwise combinations 
 \cite{DAN-2017}:
\begin{equation}
  P = \sum_{\lambda_i\in \sigma(A)} \tilde{P}_i = \sum_{\lambda_i, \lambda_j \in \sigma(A)} P_{ij} , \ \ \tilde{P}_i = \sum_{\lambda_j \in \sigma(A)} P_{ij} . \label{SD-gen-alg}
\end{equation}
If the spectrum of the dynamics matrix is simple, $\sigma(A)=\{\lambda_1,\lambda_2,\cdots,\lambda_n\}$, then 
the spectral components of SDSE of the infinite controllability Gramian is computed as  \cite{DAN-2017}
 \begin{gather}  
 \tilde{P}_i =  \left\{ R_i BB^T \, (I-\lambda_i A^T)^{-1} \right\}_H \, , \nonumber \\
  P_{ij} = \left\{ \frac{1}{1-\lambda_i \lambda^*_j} \, R_i BB^T R^*_j \right\}_H . \label{SG-gen-alg}
\end{gather}
where  
$\{\cdot\}_H$ represents the Hermitian part of a matrix, and $R_i, R_j$ are the matrix residues, that is, the coefficients in the expansion of the resolvent of matrix $A$ corresponding to its eigenvalues $\lambda_i$ and $\lambda_j$:
 \begin{equation}  \label{simple-residues}
 (Is-A)^{-1} = \sum_i \frac{R_i}{s-\lambda_i} \, .  
\end{equation}
All summation indices here and elsewhere run from $1$ to $n$, unless otherwise specified.
The spectral components $\tilde{P}_i$ and $P_{ij}$ in \eqref{SG-gen-alg} were named {\it sub-Gramians} and {\it pair sub-Gramians}, respectively. 

In this study we adopt the following basic assumptions. 
\vskip 2mm
\begin{itemize}
\item 
The Lyapunov operator in \eqref{Lyap-algeb} and \eqref{Lyap-diff} is non-singular, that is, 
all eigenvalues of matrix $A$ satisfy the condition:
\begin{equation} \label{Lyap-cond}
\lambda_i \lambda_j \ne 1 \ \ \text{for all} \ \  \lambda_i,\lambda_j\in \sigma(A) , 
\end{equation}
which is necessary and sufficient for the existence of a unique solution to these equations.
\vskip 1mm
\item 
Throughout this study, except in Section VI, the matrix $A$ is diagonalizable over the field of complex numbers and has a simple spectrum, that is, there exists an invertible matrix $U$ and a diagonal matrix $\Lambda$ with distinct numbers on the diagonal such that $A = U\Lambda U^{-1}$.
\vskip 1mm
\item 
The system \eqref{sys-0} or \eqref{sys-c} is completely controllable whenever the SDSE in companion form or the SDSE of inverse Gramians is considered.
\vskip 1mm
\item 
The matrix $A$ can be stable or unstable. For a strictly stable $A$, the solutions of \eqref{Lyap-algeb} and \eqref{Lyap-diff} can be interpreted as {\it the infinite and finite controllability Gramians}, respectively. In the second case, the obtained expansions refer to the solutions of the corresponding matrix equations.
\end{itemize}
\vskip 5mm

\section{SDSE of the solution to the difference Lyapunov equation}

The SDSE of {\it the infinite} Gramian in \eqref{SG-gen-alg} enables analysis of only the asymptotic spectral properties of the system. To study these properties over an arbitrary time interval and predict their future dynamics, 
the SDSE of {\it the finite} controllability Gramian is derived, which is a solution
to \eqref{Lyap-diff} for an arbitrary initial condition. \\
 
 \noindent \textbf{Theorem 1.} \textit{
Let system \eqref{sys-0} have a simple spectrum. Then, the SDSE of \textbf{the finite} controllability Gramian $P(t)$,
satisfying \eqref{Lyap-diff} with initial condition $P(0)=P_0$, is represented as}
\begin{gather}  
  P(t) = \sum_{i} ( \tilde{P}_i(t) + R_i P_0 (\lambda_i A^T)^{t} ) = \nonumber \\
   \sum_{i,j} ( P_{ij}(t) + R_i P_0 R^*_j (\lambda_i \lambda^*_j)^{t}), \ \  t \in \mathbb{N}, \label{SD-diff-gen} \\
  \tilde{P}_i(t) =  \left\{ R_i BB^T (I-\lambda_i A^T)^{-1} (I-(\lambda_i A^T)^{t}) \right\}_H ,   \nonumber \\
   P_{ij}(t)= \left\{ \frac{1-(\lambda_i \lambda^*_j)^{t}}{1-\lambda_i \lambda^*_j} \, R_i BB^T R^*_j \right\}_H , \label{SG-diff-gen} \\
    \tilde{P}_i(t) = \sum_{j} P_{ij}(t) . \nonumber
 \end{gather}
 
 \noindent \textbf{Proof.} If $t=0$ we obtain from \eqref{SG-diff-gen} and \eqref{SD-diff-gen} 
\begin{equation}
\tilde{P}_i (0) = 0, \     P_{ij}(0) = 0, \ 
P(0) = \sum_{i} R_i P_0 = \sum_{i,j} R_i P_0 R^*_j \, , \nonumber
\end{equation}
which gives $P(0)=P_0$ because of $\sum_{i} R_i = \sum_{j} R^*_j = I$.

Substituting $s= \lambda^{-1}_i$ into the transposed resolvent expansion \eqref{simple-residues}
$(Is-A^T)^{-1} = \sum_{j} \frac{R^*_j}{s-\lambda^*_j}$, we obtain $(I-\lambda_i A^T)^{-1} =  \sum_{j} \frac{R^*_j}{1-\lambda_i \lambda^*_j}$.
Therefore,
\begin{gather}
\tilde{P}_i (t) = \left\{ R_i BB^T (I-\lambda_i A^T)^{-1} (I-(\lambda_i A^T)^{t}) \right\}_H \nonumber \\
= \sum_{j} \left\{ \frac{R_i BB^T R^*_j}{1-\lambda_i \lambda^*_j} (I-(\lambda_i A^T)^{t}) \right\}_H \nonumber \\
= \sum_{j} \left\{ \frac{R_i BB^T R^*_j}{1-\lambda_i \lambda^*_j} (I-(\lambda_i \lambda^*_j)^{t}) \right\}_H = \sum_j P_{ij} (t) . \nonumber 
 \end{gather}
 Using the following properties of the matrix residues 
\begin{equation}
 AR_i = \lambda_i R_i, \ \ \sum_{i} R_i = \sum_{j} R^*_j = I , \nonumber
\end{equation}
 one can verify by direct substitution that the sum of terms with $P_0$ in \eqref{SD-diff-gen} satisfy \eqref{Lyap-diff} with $B=0$, 
 and the sum of terms with $BB^T$ satisfy it with $P_0 = 0$, which confirms \eqref{SD-diff-gen}.
 $\blacksquare$\\
 
 \noindent \textbf{Remark 1.} The coefficients at spectral components $\tilde{P}_i(t)$ are uniquely defined as coefficients at different geometric progressions in the solution to \eqref{Lyap-diff}, provided that the corresponding eigenvalues are distinct.
The coefficients at spectral components 
$P_{ij} (t)$ may be defined not-uniquely for some indices $i,j,k,l$ such that
$\lambda_i \lambda_j = \lambda_k \lambda_l$. However, in this case, 
the sum of such coefficients is determined uniquely.
\vskip 6mm

\section{SDSE of Gramians in Companion Form}

When expressed {\it in canonical controllability or observability form}, Gramian matrices take the simplest form, independent of the input and output matrices. In this section, we therefore derive the corresponding SDSEs for the solutions to the discrete-time algebraic and difference Lyapunov equations, as presented in Theorems~2 and~3 respectively.

\subsection{The discrete-time algebraic Lyapunov equation}

\noindent \textbf{Theorem 2.} \textit{
The SDSE of \textbf{the infinite} controllability Gramian $P_C$ in companion form in \eqref{Lyap-algeb-c} is given by}
\begin{gather}  
  P_C = \sum_{i} \tilde{P}^C_i = \sum_{i,j} P^C_{ij} , \ \ \tilde{P}^C_i = \sum_{j} P^C_{ij} ,   \label{SD-alg-c} \\
\tilde{P}^C_i = \left\{ \frac{x_i x^T_i \mathcal{I} }{\lambda^n_i N'(\lambda_i)  N(\lambda^{-1}_i)} \right\}_H \, , \label{subG-alg-c} \\
P^C_{ij} = \left\{ \frac{1}{1-\lambda_i \lambda^*_j} \, \frac{x_i x^*_j}{N'(\lambda_i) N'(\lambda^*_j)} \right\}_H \, , \label{subG2-alg-c}\\
\text{where} \ \ x_i =
   \begin{bmatrix}
   1 \\
   \lambda_i \\
   \vdots \\
   \lambda^{n-1}_i
   \end{bmatrix} \, , \
 \mathcal{I} = 
 \begin{bmatrix}
0      & \cdots    &   0     & 1  \\ 
\vdots   & \iddots   & \iddots  &0 \\
 0    & 1 & \iddots     &  \vdots \\
1  & 0 &   \cdots  &0 
\end{bmatrix} \, \nonumber 
\end{gather}
\textit{is an antidiagonal identity matrix, $N(s)=\text{det}(Is-A)$ and $N'(s)$ are the characteristic polynomial and its derivative, respectively.}\\

\noindent \textbf{Proof.} Let us find the right eigenvector $x_i$ 
for the eigenvalue $\lambda_i$ of $A_C$ in \eqref{sys-c}. 
The first component is selected as~1, and the other components are directly obtained from $A_C \, x_i= \lambda_i x_i $ 
as
\begin{equation} \label{right-ev}
x_i = [1,\lambda_i, \lambda^2_i, \cdots, \lambda^{n-1}_i]^T.
\end{equation}
For the left eigenvector $y_i$ corresponding to $\lambda_i$, we set the last component to $-1$, 
and then solve the system $y^T_i A_C = \lambda_i y^T_i$ to obtain the other components. The result is as follows: 
\begin{equation} \label{left-ev}   
y_i =\frac{1}{\lambda^n_i} \mathcal{H}_l x_i, \ \mathcal{H}_l = 
\begin{bmatrix}
0      & \cdots    &   0     & a_0  \\ 
\vdots   & \iddots   & \iddots  & a_1 \\
 0    & a_0 & \iddots     &  \vdots \\
a_0  & a_1 &   \cdots  &a_{n-1} 
\end{bmatrix} \, ,
\end{equation}
where $\mathcal{H}_l$ is a lower Hankel matrix formed by the coefficients $a_i$ of the characteristic polynomial.
The scalar product of the right and left eigenvectors is
\begin{equation} 
y^T_i x_i = - \lambda^{-1}_i \sum^n_{k_0=1} \sum^n_{k=k_0} a_k \lambda^k_i 
= -  \sum^n_{k=1} k a_{k} \lambda^{k-1}_i = - N'(\lambda_i) . \nonumber
\end{equation}
The matrix residue of the resolvent corresponding to $\lambda_i$ can be found as
\begin{equation} \label{res-0}
R_i = \frac{x_i y_i^T}{y_i^T x_i} = \frac{x_i y_i^T}{-N'(\lambda_i)} .
\end{equation}
$P_{ij}^C$ is found, according to \eqref{SG-gen-alg}. Substituting into \eqref{SG-gen-alg}:
\begin{equation}
R_i = \frac{x_i y_i^T}{-N'(\lambda_i)} , \ B =b_C, \ y^T_i b_C = -1, \ b^T_C (y^T_i)^* = -1, \ \nonumber \\
\end{equation}
obtains \eqref{subG2-alg-c} after symmetrization.
The sub-Gramian $\tilde{P}^C_i$ is found according to \eqref{SG-gen-alg}. In this formula
\begin{equation}
R_i B = \frac{x_i y^T_i b_C}{-N'(\lambda_i)} = \frac{x_i}{N'(\lambda_i)}, \ B^T = b^T_C . \label{Theo2-1}
\end{equation}
The vector ${z}_i^T=[z_1,z_2,…,z_n] = b_C^T (I - \lambda_i A_C^T )^{-1}$ we find directly from the condition
\begin{equation}
{z}_i^T (I - \lambda_i A_C^T ) = [0,\cdots,0,1] \nonumber
\end{equation}
Solving this system element by element, we obtain
\begin{equation}     
{z}_i^T = \frac{[\lambda^{-1}_i, \lambda^{-2}_i, \cdots, \lambda^{-n}_i ] }{N(\lambda^{-1}_i )} 
=\frac{x^T_i \mathcal{I}} {\lambda^n_i N(\lambda^{-1}_i )} .   \label{z_i}   
\end{equation}
Substituting \eqref{Theo2-1} and \eqref{z_i} into \eqref{SG-gen-alg} obtains \eqref{subG-alg-c}. 
$\ \blacksquare$\\

\noindent \textbf{Property 1 (Structure).} {\it The sub-Gramians $\tilde{P}^C_i$ in \eqref{SD-alg-c} are symmetric Toeplitz matrices.}\\

Indeed, it follows from \eqref{subG-alg-c} that  
$(\tilde{P}_i)_{kl} = \left\{ \frac{ \lambda_i^{k-1} \lambda_i^{n-l} }{N'(\lambda_i) \lambda^n_i N(\lambda^{-1}_i)} \right\}_H 
= \frac{ \lambda_i^{k-l} + \lambda_i^{l-k} }{2 \lambda_i N'(\lambda_i)  N(\lambda^{-1}_i)}$.\\

It is known that the infinite controllability Gramian $P_C$ for continuous LTI system in the canonical form has a special zero-plaid Hankel alternating structure \cite{Sreeram-1991}. 
In a similar way, from Property~1 we obtain that for discrete-time LTI system the controllability Gramian has a special structure. \\

\noindent \textbf{Corollary 1.} {\it The infinite controllability Gramian $P_C$ for discrete-time LTI system \eqref{Lyap-algeb-c} in the canonical form is
the symmetric Toeplitz matrix}
\begin{equation}
\tilde{P}^C = 
\begin{bmatrix}
p_1 & p_2 & \cdots & p_n \\
p_2 & p_1 & \ddots &\vdots \\
\vdots & \ddots & \ddots &p_2 \\ 
p_n & \cdots &  p_2 & p_1    
\end{bmatrix} 
\end{equation}
{\it with coefficients $p_k = \sum_{i} \frac{ \lambda_i^{k-1} + \lambda_i^{1-k} }{2 \lambda_i N'(\lambda_i)  N(\lambda^{-1}_i)}$.}\\

The pair sub-Gramians $P^C_{ij}$ in \eqref{subG-alg-c} no longer have such a Teoplitz structure, but the following property holds for them.\\

\noindent \textbf{Property 2.} {\it If system \eqref{sys-c} is stable, then $P^C_{ii} \ge 0, \, i = \overline{1,n}$.}\\

The proof follows from the fact that the absolute value of the eigenvalues of a stable system is less than one, and $P^C_{ii}$ is
proportional to $[1,\lambda_i ,\lambda_i^2,\dots, \lambda_i^{n-1}]^T ·([1,\lambda_i ,\lambda_i^2,\dots, \lambda_i^{n-1}]^T )^* \ge 0$.\\

\noindent \textbf{Example 1} illustrates an application of Theorem~2. 
Consider an unstable system \eqref{sys-c} with eigenvalues $\lambda_1=2,\lambda_2=\frac{1}{3},\lambda_3=\frac{1}{4}$
(notice that $\lambda_i \, \lambda_j \ne1 $ for all $i,j=\overline{1,3}$). 
Its characteristic polynomial and dynamics matrix take the form:
\begin{equation}
N(s)=s^3-\frac{31}{12}s^2+\frac{5}{4}s-\frac{1}{6} , \ \ A_C = 
\begin{bmatrix}
0 & 1 & 0 \\
0 & 0 & 1 \\
\frac{1}{6} & -\frac{5}{4} & \frac{31}{12}    
\end{bmatrix} . \nonumber
\end{equation}
According to \eqref{subG-alg-c}, we compute the sub-Gramians:
\begin{gather}
\tilde{P}^C_1 = -\frac{12}{35}
\begin{bmatrix}
8 & 10 & 17 \\
10 & 8 & 10 \\
17 & 10 & 8    
\end{bmatrix} , \ \tilde{P}^C_2 = -\frac{18}{55}
\begin{bmatrix}
9 & 15 & 41 \\
15 & 9 & 15 \\
41 & 15 & 9   
\end{bmatrix} , \nonumber \\
\tilde{P}^C_3 = \frac{24}{385}
\begin{bmatrix}
16 & 34 & 128.5 \\
34 & 16 & 34 \\
128.5 & 34 & 16    
\end{bmatrix} . \nonumber 
\end{gather}
Their sum gives the solution of Lyapunov equation \eqref{Lyap-algeb-c}. 
According to \eqref{subG2-alg-c}, the pair sub-Gramians are also obtained:
\begin{gather}
P^C_{11} = \frac{-48}{35^2}
\begin{bmatrix}
1 & 2 & 4 \\
2 & 4 & 8 \\
4 & 8 & 16    
\end{bmatrix} , \
P^C_{12} 
= \frac{-72}{175}
\begin{bmatrix}
18 & 21 & 37 \\
21 & 12 & 14 \\
37 & 14 & 8    
\end{bmatrix} , \ \nonumber \\
P^C_{22} = \frac{18}{25}
\begin{bmatrix}
81 & 27 & 9 \\
27 & 9 & 3 \\
9 & 3 & 1   
\end{bmatrix} , \
P^C_{23} 
= \frac{-72}{385}
\begin{bmatrix}
288 & 84 & 25 \\
84 & 24 & 7 \\
25 & 7 & 2    
\end{bmatrix} , \nonumber \\
P^C_{33} = \frac{48}{245}
\begin{bmatrix}
256 & 64 & 16 \\
64 & 16 & 4 \\
16 & 4 & 1    
\end{bmatrix} , \
P^C_{31} 
= \frac{36}{245}
\begin{bmatrix}
32 & 36 & 65 \\
36 & 16 & 18 \\
65 & 18 & 8    
\end{bmatrix} , \nonumber \\
P^C_{21} = P^C_{12}, \ \  P^C_{32} = P^C_{23}, \ \ P^C_{13} = P^C_{31}. \nonumber
\end{gather}
It can be verified that
\begin{gather}
i =\overline{1,3}: \ \tilde{P}^C_i = \sum^3_{j=1} P^C_{ij}, 
\ \ P_C = \sum^3_{i=1} \tilde{P}^C_i = \sum^3_{i,j=1} P^C_{ij} \ \square \ . \nonumber 
\end{gather}

\subsection{The difference Lyapunov equation}

\noindent \textbf{Theorem 3.} \textit{
The SDSE of \textbf{the finite} controllability Gramian $P_C(t)$ in companion form, 
satisfying \eqref{Lyap-diff-c} with initial condition $P_C (0)=P_0$, is represented as}
\begin{gather}  
  P_C(t) = \sum_{i} ( \tilde{P}^C_i(t) + R_i P_0 (\lambda_i A^T_C)^{t} ) = \nonumber \\
   \sum_{i,j} ( P^C_{ij}(t) + R_i P_0 R^*_j (\lambda_i \lambda^*_j)^{t}), \ \  t \in \mathbb{N}, \label{SD-diff-c} \\
  \tilde{P}^C_i (t) = \left\{ \frac{x_i x^T_i \mathcal{I} }{\lambda^n_i N'(\lambda_i) N(\lambda^{-1}_i)} (I-(\lambda_i A^T_C)^{t}) \right\}_H \, , \nonumber \\
P^C_{ij} (t) = \left\{ \frac{1-(\lambda_i \lambda^*_j)^{t}}{1-\lambda_i \lambda^*_j} \, \frac{x_i x^*_j}{N'(\lambda_i) N'(\lambda^*_j)} \right\}_H \, , \label{subG-diff-c} 
\end{gather}
\textit{where $R_i =\frac{x_i y^T_i}{-N'(\lambda_i)}$, $R^*_j = \frac{(y^T_j)^* x^*_j}{N'(\lambda^*_j)}$,
$x_i$ and $y_i$ are defined in \eqref{right-ev} and \eqref{left-ev}, respectively.}\\

\noindent \textbf{The proof} is analogous to Theorem~2 with the difference that instead of expressions \eqref{SG-gen-alg}
expressions \eqref{SG-diff-gen} for the difference Lyapunov equation are used. \\

\noindent \textbf{Remark 2.}
It follows from the proof of Theorem~2 that SDSE without symmetrization is also valid, and it takes form
\begin{gather}
P_C (t) = P_{C}(\infty) +\sum_i \left(R_i P_0 - \hat{P}_i^C\right) \, (\lambda_i A^T_C)^t = \nonumber \\ 
P_{C}(\infty) +\sum_{i,j}  \left( R_i P_0 R^*_j - \hat{P}_{ij}^C \right) \, (\lambda_i \lambda^*_j)^t ,  \label{SD-diff-c2} \\ 
\text{where} \ \ \ P_{C}(\infty) = \sum_i \hat{P}_i^C = \sum_{ij} \hat{P}_{ij}^C , \nonumber \\ 
\hat{P}_i^C = \frac{x_i x^T_i \mathcal{I} }{N'(\lambda_i) \lambda^n_i N(\lambda^{-1}_i)}, \nonumber \\ 
\hat{P}_{ij}^C = \frac{1}{1-\lambda_i \lambda^*_j} \, \frac{x_i x^*_j}{N'(\lambda_i) N'(\lambda^*_j)} . \label{subG-diff-c-non-sym}
\end{gather}\\

\noindent \textbf{Example 2} illustrates SDSE of the finite Gramian in companion form. 
According to Theorem 3 for the system considered in Example 1, 
the solution \eqref{Lyap-diff-c} with zero boundary condition $P_C (0)=0$ 
is decomposed over the pair spectrum in the form
\begin{gather}
P_C (t) = P_{11}^C (1 - 4^t) + (P_{12}^C+P_{21}^C )(1 - (2/3)^t ) \nonumber \\
+ \ P_{22}^C (1- (1/9)^t ) + (P_{31}^C + P_{13}^C )(1- (1/2)^t )  \nonumber \\
+ \ (P_{23}^C+P_{32}^C )(1-(1/12)^t ) + P_{33}^C (1-(1/16)^t ) , \ t=0,1,\ldots \nonumber   
\end{gather}
where the matrices $P_{ij}^C$ were found in Example 1. 
If the initial condition, $P_C (0)=P_0$ is given in \eqref{Lyap-diff-c}, 
then according to \eqref{SD-diff-c}, an additional term is added to the solution 
\begin{gather}
\tilde{P}_C (t) = P_C (t) + P_0 (t), \ \nonumber \\
P_0 (t) = \sum^3_{i,j=1} (\lambda_i \lambda_j^* )^t  \left\{R_i P_0 R_j^* \right\}_H , \nonumber
\end{gather}
where the eigenvectors $x_i, y_i$ and residue matrices $R_i, R_j$ are obtained by \eqref{right-ev}, \eqref{left-ev}, and \eqref{res-0} $\square$.

\subsection{Relationship between SDSE of Gramian in Companion Form and 
Gramian of an arbitrary MI System}

The SDSE of the controllability Gramian for an arbitrary controllable MI system \eqref{sys-0} can be obtained
from the corresponding SDSE in the canonical controllability form. \\ 

\noindent \textbf{Theorem 4.} \textit{Let a controllable MI system \eqref{sys-0} have a simple spectrum. 
The SDSE of $P(t)$, satisfying \eqref{Lyap-diff} with the initial condition $P(0)=P_0$, takes the form:}
\begin{gather}
P(t) = \sum_{i} \tilde{P}_i (t) = \sum_{i,j} P_{ij} (t) , \nonumber \\
\tilde{P}_i (t) = \mathcal{C} (\mathcal{H}_u \tilde{P}_i^C (t) \mathcal{H}_u \otimes I_m ) \mathcal{C}^*
  +\left\{ R_i P_0 (\lambda_i A^T_C)^{t} \right\}_H ,  \nonumber \\ 
P_{ij} (t) = \mathcal{C} (\mathcal{H}_u P_{ij}^C (t) \mathcal{H}_u \otimes I_m ) \mathcal{C}^*
 +\left\{ R_i P_0 R^*_j (\lambda_i \lambda^*_j)^{t})  \right\}_H , \nonumber
\end{gather}
\textit{where $\mathcal{C}$  
is the controllability matrix, $I_m$ is a unit matrix $m\times m$, $\mathcal{H}_u$ 
is given by \eqref{transfer-matrix}, and the sub-Gramians $\tilde{P}_i^C (t)$ and $P_{ij}^C (t)$ 
are obtained by \eqref{subG-diff-c} in Theorem 3.} \\ 

\noindent \textbf{The proof} follows directly from Theorem 2 and formula \eqref{MI-P-from-Pc} derived for the Gramian in \cite{Hauksdottir-2009}. \\

\section{Spectral Decomposition of Inverse of Gramian in Companion Form}

In this section, the SDSEs of {\it the inverse of infinite and finite} controllability Gramians in companion form 
are derived in Theorems~5 and~6, respectively.

\subsection{The discrete-time algebraic Lyapunov equation}

\noindent \textbf{Theorem 5.} \textit{ 
The SDSE of \textbf{the inverse of the infinite} Gramian $P_C^{-1}$ in companion form in \eqref{Lyap-algeb-c} is given by}
\begin{gather}    
P_C^{-1} = \sum_{j} \tilde{P}_j^{-C} = \sum_{j} \hat{P}_j^{-C} = \sum_{i,j} \hat{P}_{ij}^{-C}, \label{SD-invG} \\
\ \hat{P}_j^{-C} = \sum_{i} \hat{P}_{ij}^{-C},\ 
\tilde{P}_j^{-C} = \left\{ \hat{P}_j^{-C} \right\}_H, \ P_{ij}^{-C} = \left\{ \hat{P}_{ij}^{-C} \right\}_H, \ \nonumber \\
\hat{P}_j^{-C} = \frac{\lambda^n_j N(\lambda^{-1}_j)}{N'(\lambda_j)} \ \mathcal{I}  y_j  y_j^T , \ \label{SG-invG-1} \\
\hat{P}_{ij}^{-C} = \frac{(\lambda^*_i \lambda_j)^n N((\lambda^*_i)^{-1}) N(\lambda^{-1}_j)}{N'(\lambda^*_i) N'(\lambda_j)} \cdot \frac{(y^*_i)^T y^T_j}{1-\lambda^*_i \lambda_j}, \label{SG-invG-2} \\
 y_j=\frac{\mathcal{H}_l  x_j}{\lambda_j^n} , \ 
 x_j = [1,\lambda_j, \lambda^2_j, \cdots, \lambda^{n-1}_j ]^T \ , \nonumber
\end{gather}
\textit {where $\mathcal{H}_l$ is defined in \eqref{left-ev}, 
and matrices $\hat{P}_j^{-C}$ and $\hat{P}_{ij}^{-C}$ in \eqref{SG-invG-1}, \eqref{SG-invG-2} are uniquely defined by the condition} 
\begin{equation}
\forall \ i ,j: \ \ 
\hat{P}_i^{C} \hat{P}_j^{-C} = \delta_{ij} R_i , \  \hat{P}_i^{C} = \frac{x_i x^T_i \mathcal{I} }{\lambda^n_i N'(\lambda_i) N(\lambda^{-1}_i)} \label{orto-cond} 
\end{equation}
\textit {of orthogonality of eigenparts in the expansions of $P_C$ and $P_C^{-1}$ in \eqref{subG-diff-c-non-sym} and \eqref{SD-invG}, respectively.}\\

\noindent \textbf{Proof.} 
For the time-independent part of the controllability Gramian 
in \eqref{subG-diff-c-non-sym}: 
\begin{equation}
P_C = \sum_{i} \hat{P}_i^{C} = \sum_{i} \frac{x_i x^T_i \mathcal{I} }{\lambda^n_i N'(\lambda_i) N(\lambda^{-1}_i)} \nonumber
\end{equation}
the inverse matrix is constructed in the form:
\begin{equation} 
P_C^{-1}= \sum_{j} \hat{P}_j^{-C} = \sum_j \mathcal{I}  y_j  y_j^T \ \frac{ \lambda^n_j N(\lambda^{-1}_j)}{N'(\lambda_j)} \  .  \label{invG-0} 
\end{equation}
The relations $ \mathcal{I}^2 = I$ and $x^T_i y_j = -\delta_{ij} N'(\lambda_i)$ 
are used to obtain:
\begin{gather}
P_C P_C^{-1} = \sum_{i,j} \hat{P}_i^C \hat{P}_j^{-C}  
= \sum_{i,j} 
\frac{x_i x^T_i \mathcal{I} \cdot \mathcal{I}  y_j  y_j^T \, \lambda^n_j N(\lambda^{-1}_j)}{\lambda^n_i N'(\lambda_i) N(\lambda^{-1}_i)\cdot N'(\lambda_j)}
\nonumber \\
= \sum_{i,j} \frac{- x_i \delta_{ij} y_j^T  \lambda^n_j N(\lambda^{-1}_j)}{\lambda^n_i N(\lambda^{-1}_i) N'(\lambda_j)} = \sum_i R_i = I . \nonumber
\end{gather}
Carrying out the symmetrization
\begin{equation}
P_C^{-1}=\{P_C^{-1} \}_H = \sum_{j} \{\hat{P}_j^{-C} \}_H = \sum_j \tilde{P}_j^{-C} ,  \nonumber
\end{equation}
we obtain the first equality in \eqref{SD-invG}. Let us prove that 
\begin{equation}
\hat{P}_{ij}^{-C} = R^*_i \hat{P}_{j}^{-C} . \label{th4-1} 
\end{equation}
Substituting here \eqref{SD-invG}, \eqref{SG-invG-1} and $R^*_i = \frac{(y^*_i)^T x^*_i}{-N'(\lambda^*_i)}$ shows that \eqref{th4-1} holds if
\begin{equation}
\frac{-(\lambda^*_i)^n N((\lambda^*_i)^{-1})}{1-\lambda^*_i\lambda_j} = x^*_i \mathcal{I} y_j . \label{th4-2} 
\end{equation}
In addition, \eqref{z_i} yields:
\begin{gather}
\forall \lambda_i \in \sigma(A_C) : \ \ \frac{x^T_i \mathcal{I}}{(\lambda^n_i)N(\lambda^{-1}_i)} = b^T_C (I - \lambda_i A^T_C)^{-1} \nonumber \\
\Rightarrow \lambda^*_i \in \sigma(A_C) : \frac{x^*_i \mathcal{I}}{(\lambda^*_i)^n N((\lambda^*_i)^{-1})} = b^T_C (I - \lambda^*_i A^T_C)^{-1} . \nonumber 
\end{gather}
Multiplying this on the right by $y_j$ and taking into account that from (20), it follows that
\begin{gather}
 b^T_C y_j = -1, \ (I - \lambda^*_i A^T_C) y_j = (1 - \lambda^*_i\lambda_j) y_j, \nonumber \\
\Rightarrow \ \ b^T_C (I- \lambda^*_i A^T_C)^{-1} y_j = \frac{b^T_C y_j }{1-\lambda^*_i \lambda_j} = \frac{-1}{1-\lambda^*_i \lambda_j} , \nonumber
\end{gather}
\eqref{th4-2} is obtained. Therefore \eqref{th4-1} is verified. From \eqref{th4-1}, we obtain that $\hat{P}_j^{-C} = \sum_{i} \hat{P}_{ij}^{-C}$ and the validity of 
the last equality in \eqref{SD-invG}.
For matrices $\hat{P}_j^{-C}$ in \eqref{SG-invG-1} the condition \eqref{orto-cond} clearly holds. 
Let us prove the uniqueness of the representation \eqref{SG-invG-1} based on the condition \eqref{orto-cond}. 
It follows from \eqref{orto-cond} that $\hat{P}_j^{-C}$ is orthogonal to $(n-1)$ independent vectors $x^T_i \mathcal{I}, \ i\ne j$, 
that is, $\mathrm{rank} \ \hat{P}_j^{-C} = 1$. So the following representation $ \hat{P}_j^{-C} = a \, b^T $
through some non-zero vectors $a, b \in \mathbb{C}^n$ is valid.
The vector $a$ must be orthogonal to $(n-1)$ vectors $x^T_i \mathcal{I}, \ i\ne j$. 
Hence, there must be $a = \alpha \cdot \mathcal{I} y_j$ for some non-zero $\alpha \in \mathbb{C}$. There must also be 
\[
 \hat{P}_j^{C} \hat{P}_j^{-C} = 
 \frac{x_j x^T_j \mathcal{I} }{N'(\lambda_j) \lambda^n_j N(\lambda^{-1}_j)}
 \cdot \alpha \mathcal{I} y_j b^T 
=\frac{- \alpha x_j b^T} {\lambda^n_j N(\lambda^{-1}_j)} = R_j \, .
\]
Substituting here \eqref{res-0} for $R_j$, we obtain 
\[
b^T =  \frac{y^T_j \lambda^n_j N(\lambda^{-1}_j)}{\alpha N'(\lambda_j)} ,
\]
and $\hat{P}_j^{-C}$ must satisfy \eqref{SG-invG-1}. 
The uniqueness of \eqref{SG-invG-1} is proved. 
The uniqueness of \eqref{SG-invG-2} then follows from \eqref{th4-1} $\blacksquare$.\\

\noindent \textbf{Example 3} illustrates SDSE of the inverse of Gramian. 
Consider the system from Example~1. 
The inverse Gramian SDSE components are found by \eqref{SG-invG-1}:
\begin{gather}
\tilde{P}_1^{-C} = -\frac{1}{1680}
\begin{bmatrix}
24 & -91 & 145 \\
 -91 & 98 & -91 \\
 145 & -91  & 24    
\end{bmatrix} , \nonumber \\
\tilde{P}_2^{-C} = -\frac{11}{90}
\begin{bmatrix}
8 & -27 & 10 \\
-27 & 81 & -27 \\
10 & -27 & 8    
\end{bmatrix} , \ \nonumber \\
\tilde{P}_3^{-C} = \frac{55}{336}
\begin{bmatrix}
12 & -35 & 13 \\
 -35 & 98 & -35 \\
13 &  -35 & 12    
\end{bmatrix}  . \nonumber 
\end{gather} 
The pair spectral components are obtained by \eqref{SG-invG-2}:
\begin{gather}
P^{-C}_{11} = \frac{-1}{14700}
\begin{bmatrix}
1 & -7 & 12 \\
-7 & 49 & -84 \\
12 & -84 & 144    
\end{bmatrix} , \ \nonumber \\
P^{-C}_{12} = P^{-C}_{21} = \frac{11}{1050}
\begin{bmatrix}
4 & -23 & 28 \\
-23 & 126 & -136 \\
28 & -136 & 96    
\end{bmatrix} , \ \nonumber \\
P^{-C}_{22} =  \frac{121}{450}
\begin{bmatrix}
4 & -18 & 8 \\
-18 & 81 & -36 \\
8 & -36 & 16   
\end{bmatrix} , \ \nonumber \\
P^{-C}_{23} = P^{-C}_{32} = -\frac{11}{42}
\begin{bmatrix}
8 & -32 & 14 \\
-32 & 126 & -55 \\
14 & -55 & 24    
\end{bmatrix} , \nonumber \\
P^{-C}_{33} = \frac{605}{588}
\begin{bmatrix}
4 & -14 & 6 \\
-14 & 49 & -21 \\
6 & -21 & 9    
\end{bmatrix} , \ \nonumber \\
P^{-C}_{31} = P^{-C}_{13} = \frac{-11}{784}
\begin{bmatrix}
4 & -21 & 27 \\
-21 & 98 & -105 \\
27 & -105 & 72    
\end{bmatrix} . \nonumber 
\end{gather}
It can be verified that 
\begin{gather}
i =\overline{1,3}: \tilde{P}^{-C}_i = \sum^3_{j=1} P^{-C}_{ij}, 
P^{-1}_C = \sum^3_{i=1} \tilde{P}^{-C}_{i} = \sum^3_{i,j=1} P^{-C}_{ij} \ \square . \nonumber 
\end{gather}

\subsection{The discrete-time difference Lyapunov equation}

To study the spectral properties of the inverse of Gramian on an arbitrary time interval and to analyze the dynamics of the system at future times, the SDSE of the inverse of the {\it finite Gramian in companion form} is derived.\\

\noindent \textbf{Theorem 6.} \textit{Let system \eqref{sys-c} 
have a simple spectrum.
The SDSE of the inverse of \textbf{the finite} controllability Gramian 
$P_C^{-1} (t)$ in \eqref{Lyap-diff-c} with the boundary condition $P_C (0)=P_0$, takes the following form}:
\begin{gather}
P_C^{-1}(t) = G(t) \sum_j \hat{P}_j^{-C} = G(t) \sum_{ij} \hat{P}_{ij}^{-C}, \ \label{SD-invFG} \\
\text{where} \ \ G^{-1}(t) = I - \mathcal{I}(A_C^T)^t \mathcal{I} (A_C^T)^t + \sum_i \hat{P}_i^{-C}P_0 (\lambda_i A_C^T)^t \nonumber
\end{gather}
\textit{is the normalization matrix, $\hat{P}_j^{-C}$ and $\hat{P}_{ij}^{-C}$ are defined in \eqref{SG-invG-1}, \eqref{SG-invG-2},  
and $\mathcal{I}$ is an antidiagonal identity matrix.}\\

\noindent \textbf{Proof.} 
At $t=0$, it can be directly checked that $G(0) = \left( \sum_i P^{-C}_i \right)^{-1} P^{-1}_0$ and $P^{-1}_C(0) = P^{-1}_0$.
Given that $\frac{\lambda^n_j N(\lambda^{-1}_j) y_j^T x_i}{-N'(\lambda_i) \lambda^n_i N(\lambda^{-1}_i)} = \delta_{ij}$ is the Kronecker symbol, 
it is observed that
\begin{gather}
\hat{P}_j^{-C} \hat{P}_i^C = \frac{\mathcal{I} y_j y_j^T \lambda^n_j N(\lambda^{-1}_j)}{N'(\lambda_j)}
\cdot \frac{x_i x_i^T \mathcal{I}}{N'(\lambda_i) \lambda^n_i N(\lambda^{-1}_i)} \nonumber \\
= \frac{\mathcal{I} y_j \delta_{ji} x_i^T \mathcal{I}}{-N'(\lambda_j)} = \mathcal{I} R_j^T \delta_{ji} \mathcal{I}, \ \nonumber \\
\hat{P}_j^{-C} R_i = \frac{\mathcal{I} y_j y_j^T \lambda^n_j N(\lambda^{-1}_j)}{N'(\lambda_j)} \cdot \frac{x_i y_i^T}{-N'(\lambda_i)} \ \nonumber \\
= \frac{\mathcal{J} y_j N(-\lambda_j ) \delta_{ji} y^T_j}{-N'(\lambda_j)} = \delta_{ji} \hat{P}_j^{-C} \ . \nonumber
\end{gather}
The expansion for $P_C^{-1}(t)$ is checked by direct substitution:
\begin{gather}
G^{-1}(t) P_C^{-1} (t) P_C (t) = \sum_{j,i} \hat{P}_j^{-C} \hat{P}_i^C (I - (\lambda_i A_C^T)^t) \nonumber \\
+ \ \sum_{j,i} \hat{P}_j^{-C} R_i P_0 (\lambda_i A_C^T)^t \nonumber \\
=  \sum_{j,i} \mathcal{I} R_j^T \delta_{ji} \mathcal{I} (I - (\lambda_i A_C^T)^t)
+ \delta_{ji} \hat{P}_j^{-C} P_0 (\lambda_i A_C^T)^t \nonumber \\
= \mathcal{I}  \sum_{i} R_i^T \mathcal{I} - \sum_i \mathcal{I} R^T_i \lambda^t_i \mathcal{I} (A_C^T)^t \nonumber \\
+ \sum_i \hat{P}_i^{-C}P_0 (\lambda_i A_C^T)^t 
= G^{-1}(t) \ . \ \ \ \blacksquare \nonumber
\end{gather} 
\vskip3mm

\section{Extending Results to Multiple Eigenvalues} 

This section extends the previously obtained results to systems with a non-diagonalizable dynamics matrix.
This enables an SDSE of Gramians to be used to obtain closed-form estimates of the system's spectral properties 
in limiting regimes where the eigenvalues converge, such as in the case of resonant interactions.
First, we extend the SDSE of Gramians \eqref{SD-diff-gen} in Theorem~1.\\ 

\noindent \textbf{Theorem 7.} \textit{Let the eigenvalues $\lambda_i$ of matrix $A$ in \eqref{sys-0} have multiplicities $n_i$ $(n_1+\cdots+n_m =n)$.
Then SDSE of the infinite and finite Gramians $P(\infty)$ and $P(t)$ in \eqref{Lyap-algeb} and \eqref{Lyap-diff} with zero initial condition are}
\begin{gather}
P(\infty) = \sum_{i=1}^m \sum_{k=1}^{n_i} \hat{A}_k^{(i)} \cdot BB^T \cdot (A^T)^{(k-1)}(I-\lambda_i A^T)^{-k} , \nonumber \\
P(t) = \sum_{i=1}^m \tilde{P}_i (t)= \sum_{i=1}^m \sum_{k=1}^{n_i} \hat{A}_k^{(i)} \cdot BB^T \cdot  \nonumber \\
\frac{1}{(k-1)!} \,
\frac{d^{(k-1)}}{d\lambda^{(k-1)}} 
\left[(I-\lambda A^T)^{-1} (I-(\lambda A^T)^{t}) \right]_{\lambda=\lambda_i}  , \label{mult-eig-0} \\
\hat{A}_k^{(i)} = 
\frac{1}{(n_i-k)!} \cdot \frac{d^{(n_i-k)}}{d\lambda^{(n_i-k)}} \left[(\lambda-\lambda_i)^{n_i} (I\lambda-A)^{-1} \right]_{\lambda=\lambda_i} , \nonumber
\end{gather}
\textit{where $\hat{A}_k^{(i)}$ are the matrix coefficients in the partial fractional decomposition of the resolvent.}\\
 
\noindent \textbf{Proof.} The zero initial condition can be checked directly from \eqref{mult-eig-0}.
By direct substitution we can verify that
\begin{equation}
\forall i, \, k =\overline{1,n_i}: \ \ A \hat{A}_k^{(i)} =\lambda_i \hat{A}_k^{(i)} + \hat{A}_{k+1}^{(i)} (1-\delta_{k n_i}), \label{T7.1}
\end{equation}
where $\delta_{k n_i}$ is the Kronecker delta. 
Let us denote $Q = BB^T$. Using \eqref{T7.1}, we substitute $\tilde{P}_i = \sum_{k=1}^{n_i} \hat{A}_k^{(i)} Q (A^T)^{(k-1)} (I-\lambda_i A^T)^{-k} $ into \eqref{Lyap-algeb}:
\begin{gather}
A \tilde{P}_i A^T - \tilde{P}_i = \sum^{n_i}_{k=1} \hat{A}_k^{(i)} Q (I-\lambda_i A^T)^{-k} (A^T)^{(k-1)}(\lambda_i A^T - I) \nonumber \\
+ \sum^{n_i}_{k=1} \hat{A}_{k+1}^{(i)} Q (I-\lambda_i A^T)^{-k} (A^T)^{k} . \nonumber
\end{gather}
After replacing the index $\tilde{k}=k+1$ in the second sum and 
performing mutual cancellations, we obtain
\[
A \tilde{P}_i A^T - \tilde{P}_i =  - \hat{A}_1^{(i)} Q \ \Rightarrow \  A P A^T - P = - \sum_i \hat{A}_1^{(i)} Q  = -Q \, ,  
\]
where the condition $\sum_i \hat{A}_1^{(i)} = I$ follows from comparing the coefficient of $1/s$ in the large-$s$ expansion 
of the identity 
\[
(sI-A)^{-1} = \sum^m_{i=1}  \sum^{n_i}_{k=1} \frac{\hat{A}_k^{(i)}}{(s-\lambda_i)^k} \, .
\]
The expansion for $P(\infty)$ is proven.
The time-invariant parts in \eqref{mult-eig-0} coincide with the eigenparts in the expansion for $P(\infty)$ and their sum satisfies the difference equation \eqref{Lyap-diff-c}.
Let us check that the time-dependent parts of $\tilde{P}_i (t)$ satisfy the homogeneous equation 
$\tilde{P}^{(t+1)}_i = A \tilde{P}^{(t)}_i A^T$.  
Indeed, using \eqref{T7.1}, we obtain
\begin{gather}
-\sum^{n_i}_{k=1}  \frac{\hat{A}_k^{(i)} Q}{(k-1)!} \, \frac{d^{(k-1)}}{d\lambda^{(k-1)}} 
\left[(I-\lambda A^T)^{-1} (\lambda A^T)^{t+1} \right]_{\lambda=\lambda_i} = \nonumber \\
-\sum^{n_i}_{k=1}  \frac{\hat{A}_k^{(i)} Q}{(k-1)!} \, \frac{d^{(k-1)}}{d\lambda^{(k-1)}} 
\left[(I-\lambda A^T)^{-1} (\lambda A^T)^{t} \right]_{\lambda=\lambda_i}  \lambda_i A^T  \nonumber \\
-\sum^{n_i}_{\tilde{k}=2} \frac{\hat{A}_{\tilde{k}}^{(i)} Q}{(\tilde{k}-2)!} \, \frac{d^{(\tilde{k}-2)}}{d\lambda^{(\tilde{k}-2)}} 
\left[(I-\lambda A^T)^{-1} (\lambda A^T)^{t} \right]_{\lambda=\lambda_i}  A^T  . \nonumber 
\end{gather}
This follows from the differentiation formula
\[
\frac{d^{k-1}}{d\lambda^{k-1}} (f(\lambda) \cdot a\lambda) =
\frac{d^{k-1}f(\lambda)}{d\lambda^{k-1}}  \, a \lambda + (k-1) \frac{d^{k-2} f(\lambda)}{d\lambda^{k-2}} \, a
\]
for $f(\lambda) = (I-\lambda A^T)^{-1} (\lambda A^T)^{t}$ and each $k = \overline{2,n_i}$. $\blacksquare$\\

To obtain the SDSE of $P_C$ and $P^{-1}_C$ for the system in companion form, we will use the representation of 
$A$ through generalized eigenvectors. Let $m$ eigenvalues in spectrum $\sigma(A) = \{\lambda_1,\lambda_2,\dots,\lambda_m\}$ have multiplicities $n_1, n_2,\dots, n_m$, 
respectively ($n= \sum_i n_i$). Then $A$ can be represented in a Jordan canonical form
\begin{equation}
A = M J M^{-1}, \ J = J_1 \oplus J_2 \oplus \cdots \oplus J_m, \ \label{A-jordan-form} 
\end{equation}
where $J_i$ is a Jordan block corresponding to $\lambda_i$, $\oplus$ is a direct sum of matrices, 
and columns of matrix $M$ are composed of generalized right eigenvectors 
\begin{gather}
M = [x^{(1)}_1, \dots, x^{(1)}_{n_1}; x^{(2)}_1, \dots, x^{(2)}_{n_2}; \dots; x^{(m)}_1, \dots, x^{(m)}_{n_m}], \nonumber 
\end{gather}
where each eigenvalue $\lambda_i$ of multiplicity $n_i$ corresponds to the $i$-th Jordan chain of the generalized right eigenvectors 
$M_i = [x^{(i)}_1, \dots, x^{(i)}_{n_i}] $
determined sequentially from the condition
\begin{gather}
(A-\lambda_i I) x^{(i)}_1 = 0, \ (A-\lambda_i I)x^{(i)}_2 = x^{(i)}_1, \dots, \nonumber \\
 \ (A-\lambda_i I)x^{(i)}_{n_i} = x^{(i)}_{n_i-1}. \label{ev-J-chain}
\end{gather}
The eigenvectors in the $i$-th Jordan chain are not uniquely determined.
However, for system \eqref{sys-c} in companion form, they can be chosen uniquely in the simple analytical form.\\

\noindent \textbf{Lemma 1.} \textit{The $i$-th Jordan chain corresponding 
to the eigenvalue $\lambda_i$ of the system \eqref{sys-c} can be chosen as}
\begin{gather}
M_i = \left[x^{(i)}_1, \frac{d}{d\lambda} x^{(i)}_1, \dots, \frac{1}{(n_i-1)!}\frac{d^{n_i-1}}{d\lambda^{n_i-1}}x^{(i)}_{n_i} \right] \label{ev-J-chain-CF} \\
\text{where} \ \ \  x^{(i)}_1 = [ 1,\lambda_i, \lambda^2_i, \dots, \lambda^{n_i}_i]^T , \nonumber \\ 
(M_i)_{pq} =  \begin{cases} \binom{p-1}{q-1} \lambda^{p-q}_i, \ p\ge q \\ 0, \ else \end{cases}  \, . \label{ev-J-chain-CF2} 
\end{gather} 

\noindent \textbf{Proof.} 
An expression for $x^{(i)}_1$, up to a scalar factor, is obtained from the condition $(A_C - \lambda_i I)\, x^{(i)}_1 = 0$.
From the differentiation formula for each $k=\overline{1, n_i-1}$: 
\[
\frac{d^k}{d\lambda_i^k} (A_C - \lambda_i I) x^{(i)}_1 = (A_C - \lambda_i I) \frac{d^k}{d\lambda^k_i} x^{(i)}_1 - k \frac{d^{k-1}}{d\lambda^{k-1}_i} x^{(i)}_1
\]
it follows that
\[
(A_C - \lambda_i I) \frac{1}{k(k-1)!} \frac{d^k}{d\lambda^k_i} x^{(i)}_1 = \frac{1}{(k-1)!} \frac{d^{k-1}_i}{d\lambda^{k-1}} x^{(i)}_1 \, ,
\]
that is the eigenvalues in \eqref{ev-J-chain-CF} satisfy the condition of Jordan chain \eqref{ev-J-chain}.
An expression in \eqref{ev-J-chain-CF2} is obtained by differentiating $x^{(i)}_1 (\lambda_i)$. $\blacksquare$.\\

The generalized left eigenvectors are obtained  as rows in the matrix $M^{-1}$:
\begin{equation}
M^{-1} = [y^{(1)}_1, \dots, y^{(1)}_{n_1}; y^{(2)}_1, \dots, y^{(2)}_{n_2}; \dots; y^{(m)}_1, \dots, y^{(m)}_{n_m}]^T . \nonumber
\end{equation}
The normalization condition is satisfied as $x^{(i)}_k (y^{(j)}_l)^T = \delta_{ij} \delta_{kl}$.
The normalization conditions between blocks of right and left generalized eigenvectors are 
\begin{equation}
M^{(-1)}_i M_j =\delta_{ij} I_{n_i\times n_i}, \ \ M_i M^{(-1)}_i = \sum^{n_i}_{k=1} x^{(i)}_k (y^{(i)}_k)^T \, . \nonumber
\end{equation}

Now we are ready to extend SDSE in \eqref{SD-invG} in Theorem~5 to the case of multiple eigenvalues in the spectrum of $A_C$.\\

\noindent \textbf{Theorem~8.} \textit{The SDSE of discrete-time infinite Gramian $P_C(\infty)$ of system \eqref{sys-c} and its inverse $P^{-1}_C(\infty)$
are given by}
\begin{gather}
P_{C}(\infty) = \sum^m_{i=1} \hat{P}_i^C, \ \ P^{-1}_{C}(\infty) = \sum^m_{j=1} \hat{P}_j^{-C} , \nonumber \\ 
\hat{P}_i^C = M_i \mathcal{H}_i \mathcal{T}_i^T M^T_i  \mathcal {I}, \nonumber \\
\hat{P}_j^{-C} = \mathcal{I} (M^{(-1)}_j)^T \mathcal{T}^{-T}_j \mathcal{H}^{-1}_j M^{(-1)}_j, \label{SDSE-mult-c} 
\end{gather}
\begin{gather}
\text{where} \ \ \mathcal{T}_i = \mathcal{T}(\lambda_i), \nonumber \\ 
\mathcal{T}(\lambda) = 
\begin{bmatrix}
f & \frac{d f }{d\lambda} &  \frac{d^2 f }{2! \, d^2\lambda} & \cdots & \frac{d^{n_i-1} f}{(n_i-1)! \, d^{n_i-1}\lambda} \\
 0 & f & \frac{d f }{d\lambda} & \ddots & \vdots \\
 \vdots & \ddots & \ddots & \ddots & \vdots \\
 \vdots & & \ddots & f  & \frac{d f }{d\lambda} \\
 0 & \cdots & \cdots & 0 & f 
\end{bmatrix},   \label{matrix-T} \\
f (\lambda)= \Delta^{-1} , \ \  
\Delta (\lambda) = 1 + a_{n-1} \lambda + a_{n-2} \lambda^2 + \cdots + a_0 \lambda^n , \nonumber \\ 
\mathcal{H}_i =
\begin{bmatrix}
e^T_n y^{(i)}_1 & e^T_n y^{(i)}_2 & \cdots & e^T_n y^{(i)}_{n_i} \\
e^T_n y^{(i)}_2 & & \iddots  & 0 \\
\vdots & \iddots & \iddots &  \vdots \\
e^T_n y^{(i)}_{n_i} & 0 & \cdots &  0
\end{bmatrix} , \label{matrix-H} 
\end{gather}
\textit{$\mathcal{H}_i$ and $\mathcal{T}_i$ are the upper Hankel and Toeplitz matrices, respectively, 
$e_n$ is the last column of the unit matrix,
$M_i$ are given by \eqref{ev-J-chain-CF} or \eqref{ev-J-chain-CF2}, $\mathcal{I}$ is an antidiagonal identity matrix,   
and eigenparts in \eqref{SDSE-mult-c} satisfy the normalization condition $\hat{P}_i^C \hat{P}_j^{-C} =\delta_{ij} M_i M^{(-1)}_i$.}\\

\noindent \textbf{Proof.} 
According to Theorem~7, for system \eqref{sys-c} we obtain
\begin{equation}
\hat{P}_i^C =\left[ \hat{A}_1^{(i)} e_n,  \cdots, \hat{A}_{n_i}^{(i)} e_n \right]
\begin{bmatrix}
e^T_n ( I - \lambda_i A^T_C)^{-1} \\
\frac{1}{1!} \, \frac{d}{d\lambda_i} e^T_n ( I - \lambda_i A^T_C)^{-1} \\
\cdots \\
\frac{1}{(n_i-1)!} \frac{d^{n_i-1}}{d\lambda_i^{n_i-1}} e^T_n ( I - \lambda_i A^T_C)^{-1} 
\end{bmatrix} \label{prod-vect}
\end{equation}
Consider the first vector in \eqref{prod-vect}. 
Let us prove that it is equal to $M_i \mathcal{H}_i$. To do this, we prove that
\begin{equation}
\hat{A}_k^{(i)} = \sum^{n_i+1-k}_{l=1} x^{(i)}_l (y^{(i)}_{k-1+l})^T . \label{Ak-matrix}
\end{equation}
Matrices $\hat{A}_k^{(i)}$ are defined as coefficients in the partial fractional decomposition 
\begin{equation}
(I\lambda-A)^{-1} = \sum^m_{i=1} \sum^{n_i}_{k=1} \frac{\hat{A}_k^{(i)}}{(\lambda-\lambda_i)^k}. \nonumber
\end{equation}
Let us prove that \eqref{Ak-matrix} satisfies this definition, i.e.,
\begin{equation}
(I\lambda-A)^{-1} = \sum^m_{i=1} \sum^{n_i}_{k=1} \sum^{n_i+1-k}_{l=1} \frac{x^{(i)}_l (y^{(i)}_{k-1+l})^T}{(\lambda-\lambda_i)^k} . \label{eq-1}
\end{equation}
From the Jordan canonical form
\begin{equation}
A = M J M^{-1} = \sum^m_{i=1} \left( \sum^{n_i}_{k=1} \lambda_i x^{(i)}_k (y^{(i)}_k)^T +\sum^{n_i-1}_{k=1} x^{(i)}_k (y^{(i)}_{k+1})^T \right) \nonumber
\end{equation}
it follows that
\begin{gather}
I \lambda-A = \lambda M M^{-1} - M J M^{-1} \nonumber \\
= \sum^m_{i=1} \left( \sum^{n_i}_{k=1} (\lambda-\lambda_i) x^{(i)}_k (y^{(i)}_k)^T - \sum^{n_i-1}_{k=1} x^{(i)}_k (y^{(i)}_{k+1})^T \right) . \label{eq-2}
\end{gather}
Multiply \eqref{eq-1} by \eqref{eq-2} using the orthogonality of eigenvectors 
\begin{gather}
 \sum^m_{i,i'=1} \sum^{n_i}_{k=1} \sum^{n_i+1-k}_{l=1} \frac{x^{(i)}_l (y^{(i)}_{k-1+l})^T}{(\lambda-\lambda_i)^k} \times \nonumber \\
 \times \left[  \sum^{n_i}_{k'=1} (\lambda-\lambda_{i'}) x^{(i')}_{k'} (y^{(i')}_{k'})^T - \sum^{n_i-1}_{k'=1} x^{(i')}_{k'} (y^{(i')}_{k'+1})^T \right] \nonumber \\
 =  \sum_{i} \left[ \sum^{n_i}_{k=1} \sum^{n_i+1-k}_{l=1} \frac{x^{(i)}_l (y^{(i)}_{k-1+l})^T}{(\lambda-\lambda_i)^{k-1}}
 - \sum^{n_i-1}_{k=1} \sum^{n_i-k}_{l=1} \frac{x^{(i)}_l (y^{(i)}_{k+l})^T}{(\lambda-\lambda_i)^{k}} \right] \nonumber \\
= \sum_i \left[  
 \sum^{n_i}_{l=1} x^{(i)}_l (y^{(i)}_l)^T 
 + \sum^{n_i - 1}_{\tilde{k}=1} \sum^{n_i-\tilde{k}}_{l=1} \frac{x^{(i)}_l (y^{(i)}_{\tilde{k}+l})^T}{(\lambda-\lambda_i)^{\tilde{k}}} \right. \nonumber \\
 \left. - \sum^{n_i-1}_{k=1} \sum^{n_i - k}_{l=1} \frac{x^{(i)}_l (y^{(i)}_{k+l})^T}{(\lambda-\lambda_i)^k} \right] 
 =  \sum_i \sum^{n_i}_{l=1} x^{(i)}_l (y^{(i)}_l)^T = I . \nonumber 
\end{gather}
The obtained identity proves \eqref{Ak-matrix}; therefore, we obtain
\begin{equation}
\left[ \hat{A}_1^{(i)} e_n,  \cdots, \hat{A}_{n_i}^{(i)} e_n \right] = M_i \mathcal{H}_i . \label{1-vector}
\end{equation}

Consider the second vector $X$ in \eqref{prod-vect}.
By directly solving the equation $z^T (I-\lambda_i A^T_C) = e^T_n$ we obtain $z^T = \Delta^{-1}(\lambda_i) [\lambda_i^{n-1}, \lambda_i^{n-2}, \cdots, \lambda_i, 1]$, or in other terms $ e^T_n (I-\lambda_i A^T_C)^{-1} = \Delta^{-1}(\lambda_i) \, x^{(i)}_1 \mathcal{I}$. Substituting this to \eqref{prod-vect}, we obtain
\begin{equation} 
X = 
\begin{bmatrix}
\Delta^{-1}(\lambda_i) \, x^{(i)}_1 \\
\frac{1}{1!} \, \frac{d}{d\lambda_i} \Delta^{-1}(\lambda_i) \, x^{(i)}_1 \\
\cdots \\
\frac{1}{(n_i-1)!} \frac{d^{n_i-1}}{d\lambda_i^{n_i-1}} \Delta^{-1}(\lambda_i) \, x^{(i)}_1 
\end{bmatrix} \mathcal{I} \, . \label{sec-vector} 
\end{equation}
For each $ k = \overline{0, n_i-1}$ it holds
\[
\frac{1}{k!} \frac{d^k}{\lambda_i^k} (\Delta^{-1}(\lambda_i) x^{(i)}_1)
= \sum^k_{j=0} \frac{1}{(k-j)!} \frac{d^{k-j}}{d\lambda_i^{k-j}} (\Delta^{-1})
\frac{1}{j!} \frac{d^j}{\lambda_i^j} (x^{(i)}_1) \, .
\]
Taking this and \eqref{ev-J-chain-CF} into account, \eqref{sec-vector} can be written as  
\begin{equation}
\mathcal{I} X^T_i = M_i \mathcal{T}_i \, ,\ \ \text{or} \ \ X _i = \mathcal{T}^T_i M^T_i \mathcal{I} \, . \label{2-vector} 
\end{equation}
From \eqref{1-vector} and \eqref{2-vector} we obtain $\hat{P}_i^C$.
The expression for $\hat{P}_j^{-C}$ is chosen to satisfy the normalization condition. 
In this case
\begin{gather}
P_C P^{-1}_C =\sum^m_{i,j=1} \hat{P}_i^C \hat{P}_j^{-C} = \sum^m_{i=1} M_i M^{(-1)}_i 
= M M^{-1} = I \,  \blacksquare .  \nonumber
\end{gather}

\noindent \textbf{Example~4} illustrates and verifies the formulas of Theorem~8.
Consider an unstable system with
$\lambda_1=3$ of multiplicity 2 and $\lambda_2=1/2$ of multiplicity 2 in companion form. We have 
\begin{gather}
N(s) = s^4 - 7s^3 + 15.25 s^2 - 10.5s + 2.25 \, , \nonumber \\
a^T =[2.25,-10.5,15.25,-7] . \nonumber 
\end{gather}
According to \eqref{ev-J-chain-CF} we obtain the generalized eigenvectors
\begin{gather}
M=[M_1,M_2] =
\begin{bmatrix}
1 & 0 & 1  & 0  \\
3 & 1 & 1/2  & 1  \\
9 & 6 & 1/4  & 1  \\
27 & 27 & 1/8  & 3/4  
\end{bmatrix}, \nonumber \\
M^{-1}=
\begin{bmatrix}
M^{-1}_1 \\
M^{-1}_2 
\end{bmatrix} = \frac{1}{125}
\begin{bmatrix}
17 & -72 & 84  & -16 \\
-15 & 65 & -80  & 20 \\
108 & 72 & -84 & 16  \\
-90 & 240 & -130  & 20  
\end{bmatrix} . \nonumber
\end{gather}
According to \eqref{matrix-T} and \eqref{matrix-H} we obtain
\begin{gather}
\mathcal{T}_1 = \frac{1}{64}
\begin{bmatrix}
4 & -11 \\
0 &  4 
\end{bmatrix}, \ \mathcal{H}_1 = \frac{4}{125}
\begin{bmatrix}
-4 & 5 \\
5 & 0 
\end{bmatrix}, \nonumber \\
\mathcal{T}_2 = \frac{64}{27}
\begin{bmatrix}
3 & -32  \\
0 & 3  
\end{bmatrix}, \ \mathcal{H}_2 = \frac{4}{125}
\begin{bmatrix}
4 & 5  \\
5 & 0 
\end{bmatrix}. \nonumber 
\end{gather}
According to \eqref{SDSE-mult-c} we obtain the sub-Gramians
\begin{gather}
\hat{P}^C_1 =
 \frac{-1}{2000}
\begin{bmatrix}
1377 &  519 & 193 &  71 \\
3591 & 1377 & 519 & 193 \\
9153 &  3591  & 1377 & 519 \\
22599 & 9153 & 3591  & 1377 
\end{bmatrix}, \  \nonumber \\
\hat{P}^C_2 = \frac{-16}{3375}
\begin{bmatrix}
116 &  352 & 944 & 2368 \\
28  & 116 & 352   & 944  \\
-1   &  28  & 116 & 352 \\
-8    & -1    & 28   & 116 
\end{bmatrix}. \  \nonumber
\end{gather}
Their sum 
$\hat{P}^C_1+\hat{P}^C_2$
provides an exact solution of Lyapunov equation \eqref{Lyap-algeb-c}.
According to \eqref{SDSE-mult-c} we obtain the eigenparts for the inverse of the Gramian
\begin{gather}
\tilde{P}^{-C}_1 = \frac{1}{125} 
\begin{bmatrix}
-388 &  1708 & -2176 &  624 \\
1312 & -5792 & 7424 & -2176 \\
-1021 &  4511 & -5792 & 1708 \\
231 & -1021 & 1312 & -388 
\end{bmatrix}, \  \nonumber \\
\tilde{P}^{-C}_2 = \frac{3}{2000}
\begin{bmatrix}
-639 &  1974 & -1103 & 172 \\
4086 & -12651 & 7072   & -1103 \\
-7263 &  22608 & -12651 & 1974 \\
2268  & -7263 & 4086 & -639 
\end{bmatrix}. \nonumber
\end{gather}
Their sum 
$\hat{P}^{-C}_1+\hat{P}^{-C}_2$
provides 
$P^{-1}_C$ 
in closed-form. $\square$ \\

\section{Conclusion}
This paper developed closed-form decompositions of the discrete-time Gramian and of its inverse in terms of the eigenvalues of the dynamics matrix. The same viewpoint was then extended to the solutions of the difference Lyapunov equation and generalized to the case of multiple eigenvalues. Taken together, these results provide a mode-resolved description of discrete-time Gramian structure that complements the more familiar use of aggregate spectral quantities such as eigenvalues, singular values, determinants, and traces. 

Several directions follow naturally from these results. One is the development of Gramian-based algorithms that use modal decomposition directly rather than only through determinant, trace, or singular-value summaries. Another is the use of the proposed formulas in large-scale settings, where structure-exploiting approximations may reduce computational burden in actuator and sensor design, model reduction, or energy-aware control. It is also of interest to investigate how these decompositions interact with stochastic estimation, finite-horizon formulations, and networked dynamics, where discrete-time Gramian methods are already well established \cite{Bag-2009}, \cite{Pasq-2013}, \cite{Summ-2014}, \cite{Yan-2012}, \cite{Bag-2021}. These directions suggest that closed-form spectral decompositions of discrete-time Gramians may serve not only as analytical tools, but also as building blocks for improved algorithms in a wider class of control and estimation problems.

\section*{Acknowledgment}
We thank I.B.~Yadykin for the fruitful discussions and A.A.~Ulyutichev for verifying the examples.

\end{document}